\font \sevenrm=cmr7
\font \eightrm=cmr8

\font \eightbf=cmbx8
\font \bigbf=cmbx10 scaled \magstep1
\font \Bigbf=cmbx10 scaled \magstep2

\font \tengoth=eufm10
\font \sevengoth=eufm7
\font \fivegoth=eufm5

\newfam\gothfam
\textfont \gothfam=\tengoth
\scriptfont \gothfam=\sevengoth
\scriptscriptfont \gothfam=\fivegoth

%

\font \tenmath=msbm10
\font \sevenmath=msbm7
\font \fivemath=msbm5

\newfam\mathfam
\textfont \mathfam=\tenmath
\scriptfont \mathfam=\sevenmath
\scriptscriptfont \mathfam=\fivemath
\def\math{\fam\mathfam\tenmath}

%
%
%
\def\titre#1{\centerline{\Bigbf #1}\nobreak\nobreak\vglue 10mm\nobreak}

\def\paragraphe#1{\bigskip\goodbreak {\bigbf #1}\nobreak\vglue 12pt\nobreak}
\def\alinea#1{\medskip\allowbreak{\bf#1}\nobreak\vglue 9pt\nobreak}
\def\ssq{\smallskip\qquad}

%
%
\def\th#1{\bigskip\goodbreak {\bf Theorem #1.} \par\nobreak \sl }
\def\prop#1{\bigskip\goodbreak {\bf Proposition #1.} \par\nobreak \sl }
\def\lemme#1{\bigskip\goodbreak {\bf Lemma #1.} \par\nobreak \sl }

\def\dem{\bigskip\goodbreak \it Proof. \rm}
\def\ndem{\bigskip\goodbreak \rm}
\def\qed{\par\nobreak\hfill $\bullet$ \par\goodbreak}
%
%
\def\uple#1#2{#1_1,\ldots ,{#1}_{#2}}
\def\corde#1#2-#3{{#1}_{#2},\ldots ,{#1}_{#3}}

\def\ordcorde#1#2-#3{{#1}_{#2} \le \cdots \le {#1}_{#3}}
\def\strictordcorde#1#2-#3{{#1}_{#2} < \cdots < {#1}_{#3}}
\def \restr#1{\mathstrut_{\textstyle |}\raise-6pt\hbox{$\scriptstyle #1$}}
\def \srestr#1{\mathstrut_{\scriptstyle |}\hbox to -1.5pt{}\raise-4pt\hbox{$\scriptscriptstyle #1$}}
\def \inver{^{-1}}
\def\dbar{d\!\!\hbox to 4.5pt{\hfill\vrule height 5.5pt depth -5.3pt
        width 3.5pt}}

\def\frac#1#2{{\textstyle {#1\over #2}}}

\def\R{{\math R}}

\def\N{{\math N}}

\def\permuc#1#2#3#4{#1#2#3#4+#1#3#4#2+#1#4#2#3}

\def\fleche#1{\mathop{\hbox to #1 mm{\rightarrowfill}}\limits}
\def\gfleche#1{\mathop{\hbox to #1 mm{\leftarrowfill}}\limits}
\def\inj#1{\mathop{\hbox to #1 mm{$\lhook\joinrel$\rightarrowfill}}\limits}
\def\ginj#1{\mathop{\hbox to #1 mm{\leftarrowfill$\joinrel\rhook$}}\limits}
\def\surj#1{\mathop{\hbox to #1 mm{\rightarrowfill\hskip 2pt\llap{$\rightarrow$}}}\limits}
\def\gsurj#1{\mathop{\hbox to #1 mm{\rlap{$\leftarrow$}\hskip 2pt \leftarrowfill}}\limits}
%
%
\def \g#1{\hbox{\tengoth #1}}
\def \sg#1{\hbox{\sevengoth #1}}

\def\Cal #1{{\cal #1}}
%
%

\def \mop#1{\mathop{\hbox{\rm #1}}\nolimits}

%
%
\def \bib #1{\null\medskip \strut\llap{[#1]\quad}}
\def\cite#1{[#1]}
\magnification=\magstep1
\parindent=0cm
\def\titre#1{\centerline{\Bigbf #1}\vskip 16pt}
\def\paragraphe#1{\bigskip {\bigbf #1}\vskip 12pt}
\def\alinea#1{\medskip{\bf #1}\vskip 6pt}
\def\ssq{\smallskip\qquad}

\let\bord=\partial
\let\wt=\widetilde
\def\espace #1{\hbox to #1 pt{}}
\def\esp{\espace{15}}
\titre{On quantization of quadratic Poisson structures}
\centerline{D. Manchon
\footnote{$^1$}{\eightrm Institut Elie Cartan, CNRS, BP 239, 54506 Vandoeuvre CEDEX. manchon@iecn.u-nancy.fr}, 
M. Masmoudi
\footnote{$^2$}{\eightrm Universit\'e de Metz, D\'epartement de Math\'ematiques, \^ile du Saulcy, 57045 Metz CEDEX 01. masmou\-di@poncelet.univ-metz.fr},
A. Roux
\footnote{$^3$}{\eightrm Universit\'e de Metz, D\'epartement de Math\'ematiques, \^ile du Saulcy, 57045 Metz CEDEX 01. roux@pon\-celet.univ-metz.fr}
}
\vskip 15mm
{\baselineskip=10pt \eightbf Abstract : \eightrm Any classical r-matrix on the Lie algebra of linear operators on a real vector space $V$ gives rise to a quadratic Poisson structure on $V$ which admits a deformation quantization stemming from the construction of V. Drinfel'd \cite {Dr}, \cite{Gr}. We exhibit in this article an example of quadratic Poisson structure which does not arise this way.\par}
\vskip 6mm
\paragraphe{I. Introduction}
let $V$ be a finite-dimensional real vector space. The linear action of the Lie group $Gl(V)$ on $V$ induces by differentiation a Lie algebra isomorphism from $\g g=\g gl(V)$ to the Lie algebra of linear vector fields on $V$. Given a basis $(\uple en)$ and then identifying $\g gl(V)$ with the Lie algebra of real $n\times n$ matrices the isomorphism is given by~:
$$J(E_{ij})=x_i\bord_j,$$
where $E_{ij}$ is the matrix with entries all vanishing except one equal to $1$ on the i$^{\hbox{\sevenrm th}}$ line and j$^{\hbox{\sevenrm th}}$ column.
\ssq
There is a unique way to extend the Lie bracket of $\g g$ to a graded Lie bracket, called the {\sl Schouten bracket\/} on the shifted exterior algebra $\Lambda(\g g)[1]$ in a way compatible with the exterior product. The shift means that elements of $\Lambda^k (\g g)$ are of degree $k-1$, and then the Schouten bracket maps $\Lambda^k (\g g)\times \Lambda^l (\g g)$ to $\Lambda^{k+l-1}(\g g)$. The exterior algebra $\Lambda (\g g)$ inherits then a structure of Gerstenhaber algebra (cf. for example \cite{V}, introduction).
\ssq
The space $T^{\hbox{\sevenrm poly}}(V)$ of polyvector fields on $M$ is also endowed with a Gerstenhaber algebra structure, with Schouten bracket extending Lie bracket of vector fields \cite{V}. The subalgebra (for exterior product) generated by linear vector fields is a Gerstenhaber subalgebra $\wt\Lambda(V)$ of $T^{\hbox{\sevenrm poly}}(V)$. The isomorphism $J$ extends to a surjective Gerstenhaber algebra morphism~:
$$J^\bullet:\Lambda^\bullet(\g g)\longmapsto \wt\Lambda^\bullet(V).$$
Map $J^k$ has nontrivial kernel for $k\ge 2$ as long as $V$ has dimension $\ge 2$~: for example we have~:
$$J^2(E_{ij}\wedge E_{kj})=0.$$
A {\sl classical r-matrix\/} on $\g g$ is by definition an element $r$ of $\g g\wedge \g g$ such that $[r,r]=0$. According to the discussion above the bivector field $J^2(r)$ defines then a quadratic Poisson structure on $V$. A natural question arises then : can one recover this way any quadratic Poisson structure on $V$? It was claimed true in \cite{BR} but Z.H. Liu and P. Xu discovered that the authors' argument was not correct \cite {LX}. They brought up a positive answer in the two-dimensional case (\cite {LX} prop. 2.1) : namely the general quadratic Poisson structure $(ax_1^2+2bx_1x_2+cx_2^2)\bord_1\wedge\bord_2$ is equal to $J^{2}(r)$ with~:
$$r=\pmatrix{b & -a \cr
                c& -b\cr}\wedge\pmatrix{1&0\cr 0&1\cr}.$$
We give here a {\sl negative\/} answer to this question in general : after a somewhat lengthy but elementary computation we show in \S\ III that bivector field $(x_1^2+x_2x_3)\bord_2\wedge\bord_3$ on $\R^3$ is a counterexample to this conjecture~: it is outside the image of the set of r-matrices by $J^2$.
\ssq
We recall in \S \ II the construction by V.G. Drinfel'd of a translation-invariant deformation quantization on any Lie group $G$ once given a classical r-matrix on the Lie algebra $\g g$ \cite{Dr}, \cite{T}. The problem reduces to the case when $r$ is non-degenerate, and the deformation quantization is then obtained by suitable restriction and transportation of Baker-Campbell-Hausdorff deformation quantization (\cite {Ka}) of the dual $\wt{\g g}^*$ of the central extension $\wt{\g g}$ of $\g g$ defined by $r$. The construction works moreover for any Kontsevich-type star product \cite {ABM} on $\wt{\g g}^*$. For $\g g=\g gl(V)$, such a star product on $\wt{\g g}^*$ gives almost immediately through this construction a deformation quantization of quadratic Poisson structure $J^2(r)$ on $V$.
\ssq
Deformation quantization of some particular quadratic Poisson structures has been considered by several people namely Omori, Maeda and Yoshioka \cite{OMY prop. 4.7}. Explicit computations for all quadratic Poisson structures in dimension $3$ (then including our counterexample as well) have been performed by El Galiou and Tihami \cite {ET}, by a case-by-case method based on the classification of Dufour and Haraki \cite {DH}. Let us recall that the existence of a deformation quantization for any Poisson manifold is a direct consequence of M.Kontsevich's formality theorem.
\smallskip
{\bf Acknowledgements\/}~: the authors thank Maxim Kontsevich and Daniel Sternheimer for useful comments and precisions.
\paragraphe{II. Quantization of Poisson structures coming from r-matrices}
\qquad Let $\g g$ be a Lie algebra, and let $r\in \g g\wedge\g g$ a classical r-matrix. It defines an antisymmetric operator ~:
$$\wt r:\g g^*\longrightarrow \g g.$$ 
Classical Yang-Baxter equation $[r,r]=0$ is equivalent to~:
\def\tata#1#2#3{<#1,\,[\wt r(#2),\,\wt r(#3)]>}
$$\permuc\tata \xi\eta\zeta =0\eqno{(*)}$$
for any $\xi,\eta,\zeta\in \g g^*$. The r-matrix $r$ defines a left translation-invariant Poisson structure on any Lie group $G$ with Lie algebra $\g g$.
\alinea{II.1. A central extension}
\qquad We can firstly suppose $r$ nondegenerate, i.e. that $\wt r$ is inversible with inverse $\wt \omega$, where $\omega$ belongs to $\g g^*\wedge \g g^*$.  Classical Yang-Baxter equation is in this case equivalent to~:
\def\toto#1#2#3{<\wt\omega #1,\,[#2,#3]>}
$$\permuc\toto XYZ=0$$
for any $X,Y,Z\in \g g$, i.e is equivalent to the fact that $\omega$ is a 2-cocycle with values in the trivial representation. Let $\wt {\g g}$ the central extension of $\g g$ by this cocycle, defined by $\wt {\g g}=\g g\oplus\R$ with bracket~:
$$[X+\alpha,Y+\beta]=[X,Y]+<\omega, X\wedge Y>.$$
The cocycle condition on $\omega$ is equivalent to de Jacobi identity for this bracket. Let $X_0=(0,1)\in\wt{\g g}$, and let $\Cal H$ the hyperplane in $\wt {\g g}^*$ defined by~:
$$\Cal H=\{\xi \in \wt{\g g}^*/<\xi,X_0>=1\}.$$
It is the symplectic leaf through the point $\xi_0$ defined by $<\xi_0,\g g>=0$ and $<\xi_0,X_0>=1$. It is then the coadjoint orbit $\mop{Ad}^*\wt G.\xi_0$ for any Lie group $\wt G$ with Lie algebra $\wt{\g g}$.
\alinea{II.2. Kontsevich star products \rm(after \cite {ABM})}
The linear Poisson manifold $\wt {\g g}^*$ admits a whole bunch of equivalent $\mop{Ad}^*\wt G$-invariant deformation quantizations which can be built from the enveloping algebra $\Cal U(\wt{\g g})$, for example the Baker-Campbell-Hausdorff quantization or the Kontsevich quantization \cite{K}, \cite{ABM}, \cite{Ka}, \cite{Di}. The baker-Campbell-Hausdorff quantization is given by the following integral formula \cite {ABM}~:
$$(u\mathop\#\limits_{BCH} v)(\xi)=\int\!\!\!\int_{\sg g\times\sg g}
        \Cal F\inver u(x)\Cal F\inver v(y)e^{i<\xi,\,x\mathop .\limits_{\hbar} y>}
                \,dx\,dy,$$
where the inverse Fourier transform is given by~:
$$\Cal F\inver u(x)=(2\pi)^{-n}\int_{\sg g^*}u(\eta)e^{i<x,\eta>}d\eta,$$
and $x\mathop .\limits_{\hbar} y$ stands for the Baker-Campbell-Hausdorff expansion~:
$$x+y+{\hbar\over 2}[x,y]+{\hbar^2\over 12}([x,[x,y]]+[y,[y,x]])+\cdots.$$

The Lebesgue measure $d\eta$ on $\g g^*$ is normalized so that it is the dual mesure of Lebesgue measure $dx$ on $\g g$. The quantizations we can consider here are the ones called ``Kontsevich star products'' in \cite {ABM}. They are all equivalent to the BCH quantization. The equivalence is a formal series of differential operators with constant coefficients on $\g g^*$ precisely given by a formal series of $G$-invariant polynomials on $\g g$ of the following form~:
$$F(x)=1+\sum_{k\ge 1}\hbar^{2k}\sum_{c\ge 1}\sum_{(\uple sc)\in S_{2k}^c}a_{\uple sc}\mop{Tr}(\mop{ad}x)^{s_1}\cdots\mop{Tr}(\mop{ad}x)^{s_c},$$
where $S_{2k}^c$ stands for those $(s_1,\ldots,s_c)$ in $\N^c$ such that $s_1+\cdots +s_c=2k$, $s_1\le s_2\le\cdots\le s_c$ and $s_j\not =1$. The star product obtained this way admits the following integral form~:
$$(u\#v)(\xi)=\int\!\!\!\int_{\sg g\times\sg g}
        \Cal F\inver u(x)\Cal F\inver v(y){F(-ix)F(-iy)\over
        F\bigl(-i(x\mathop .\limits_{\hbar}y)\bigr)}e^{i<\xi,\,x\mathop .\limits_{\hbar} y>}
                \,dx\,dy.$$

\alinea{II.3. Quantization of left-invariant Poisson structures}
It is easy to derive from the fact that $X_0$ is central that any of the deformation quantizations defined above does define by restriction a deformation quantization of $\Cal H$. Let $G$ be the subgroup of $\wt G$ with Lie algebra $\g g$. We clearly have~:
$$\mop{Ad}^*G.\xi_0=\mop{Ad}^*\wt G.\xi_0=\Cal H.$$
It is moreover easy to check that the stabilizer of $\xi_0$ in $\wt G$ is the one-dimensional subgroup with Lie algebra generated by $X_0$. It is a simple consequence of the nondegeneracy of the alternate bilinear form $\omega$. The dimension of $G$ is the equal to the dimension of $\Cal H$. The map~:
$$\eqalign{\varphi:G'   &\longrightarrow \Cal H \cr
                     g  &\longmapsto \mop{Ad}^*g.\xi_0  \cr}$$
is then a local $G$-equivariant diffeomorphism near the identity (with left translation on the left-hand side and coadjoint action on the right-hand side). We can then transport any deformation quantization of $\Cal H$ and get a left translation-invariant deformation quantization of a neighbourhood of the identity in $G$. It extends by translation invariance to the whole group $G$, as well as to any Lie group $G'$ locally isomorphic to $G$.
\ssq
The deformation quantization on $G$ can be written~:
$$u\# v=\sum_{k\ge 0}\hbar^kC_k(u,v),$$
where the $C_k$'s are left-invariant bidifferential operators on $G$. There exists then an element $F=\sum \hbar^kF_k$ in $\bigl(\Cal U(\g g)\otimes\Cal U(\g g)\bigr)[[\hbar]]$ such that~:
$$u\#v(g)=F(u\otimes v)(g,g).\eqno{(**)}$$
Let us now fix a basis $\uple xn$ of $\g g$, and consider elements of $\Cal U(\g g)$ as polynomials $F(x)=F(\uple xn)$ of the $n$ noncommuting variables $\uple xn$, which satisfy the relations~:
$$x_ix_j-x_jx_i=[x_i,x_j]=\sum_kc_{ij}^kx_k.$$
Introducing a second identical set of noncommuting variables $y=(\uple yn)$ commuting with th $x_j's$ we can write any element $A\in\Cal U(\g g)\otimes\Cal U(\g g)$ as $A(x,y)$. The element $F$ defined above can then be written $F(x,y)$ as a formal series with coefficients $F_k(x,y)$.
\prop{II.1}
The formal series $F=F(x,y)\in\bigl(\Cal U(\g g)\otimes\Cal U(\g g)\bigr)[[\hbar]]$ above verifies~:
$$\eqalign{&\hbox{1) } F_0(x,y)=0.\cr
        &\hbox{2) } F_1(x,y)={1\over 2}\sum_{i,j}r_{ij}x_iy_j.\cr
        &\hbox{3) } F_k(x,0)=F_k(0,y)=0 \hbox{ for }k\ge 1.\cr
        &\hbox{4) } F(x+y,z)F(x,y)=F(x,y+z)F(y,z).\cr}$$
Conversely any $F(x,y)$ endowed with those 4 properties defines by formula $(**)$ a left translation deformation quantization of $G$.
\dem
It is well-known : see for example \cite {Dr}, \cite T : first condition comes from the fact that $C_0(u,v)$ is the ordinary product $uv$. Second property comes from the expression of left-invariant Poisson bracket on $G$ defined from the r-matrix, third property expresses the fact that $1\#u=u\#1=u$, and last property is an expression of the associativity of star product $\#$. Let us elaborate a bit on that last point : any element $X_j$ of the basis corresponds to polynomial expression $G(x)=x_j$. The Leibniz rule~:
$$X_j.(\varphi\psi)=(X_j.\varphi)\psi + \varphi.(X_j.\psi)$$
can be written as~:
$$G(x)\circ m=m\circ G(x+y),$$
where $m:C^\infty(G\times G)\rightarrow C^\infty(G)$ stands for multiplication, here the restriction to the diagonal. The formula above extends to any polynomial expression $G$ representing any element of the enveloping algebra. We have then~:
$$\eqalign{(u\#v)\#w    &=\bigl(m\circ F(u\otimes v) \bigr)\#w  \cr
        &=m\circ F\Bigl(\bigl(m\circ F(u\otimes v)  \bigr) \otimes w\Bigr) \cr
        &=m\circ F\circ (m\otimes I)\circ (F\otimes I)(u\otimes v\otimes w)\cr
        &=m\circ F(x,z)\circ(m\otimes I)\circ F(x,y)(u\otimes v\otimes w)\cr
        &=m\circ (m\otimes I)\circ F(x+y,z)F(x,y)(u\otimes v\otimes w).\cr}$$
Similarly we have~:
$$u\#(v\#w)=m\circ (I\otimes m)\circ F(x,y+z)F(y,z)(u\otimes v\otimes w).$$
The associativity condition for product $\#$ is then equivalent to property 4) of the proposition.
\qed 
let us now look at the case when $r$ is degenerate. Then the image $\g g_0$ of $\wt r$ is a subspace strictly contained in $\g g$. By skew-symmetry $\g g_0$ is also the orthogonal of the kernel of $\wt r$, and classical Yang-Baxter equation $[r,r]=0$ ensures thanks to $(*)$ that $\g g_0$ is a Lie subalgebra of $\g g$. We get this way a nondegenerate $r_0\in\g g_0\wedge\g g_0$ such that $[r_0,r_0]=0$. Applying the procedure above we get an $F=\sum \hbar^kF_k$ in $\bigl(\Cal U(\g g_0)\otimes\Cal U(\g g_0)\bigr)[[\hbar]]$ which can be seen as an element of
$\bigl(\Cal U(\g g)\otimes\Cal U(\g g)\bigr)[[\hbar]]$.
\alinea{II.4. A class of easily quantizable Poisson structures}
\qquad Let $G$ be a Lie group with Lie algebra $\g g$. Let $F(x,y)$ a formal series in $\Cal U(\g g)\otimes\Cal U(\g g)[[\hbar]]$ satisfying properties 1-4 of proposition II.1 (for example that one constructed from an r-matrix along the lines above). Let $M$ be any differentiable manifold endowed with an action of $G$. The differentiation of this action induces a Lie algebra morphism from $\g g$ to the vector fields on $M$, which extends to an algebra morphism from $\Cal U(\g g)$ to the algebra of differential operators on $M$. Similarly it induces an algebra morphism from $\Cal U(\g g)\otimes\Cal U(\g g)$ to the algebra of differential operators on $M\times M$. The formal series of bidifferential operators defined by the formula~:
$$*=m\circ F(x,y)$$
(where $m:C^\infty(M\times M)\rightarrow C^\infty(M)$ stands for ordinary multiplication of functions on $M$) defines then a star product on $M$, the associated Poisson bivector being defined by $F_1(x,y)-F_1(y,x)$. The proof of this fact goes the same way as that of proposition II.1. It is easily seen that if $F(x,y)$ comes from a classical r-matrix $r\in \g g\wedge\g g$ then the Poisson structure on $M$ is $J^2(r)$ where $J^\bullet$ is the Gerstenhaber algebra morphism from $\Lambda (\g g)$ to multivector fields on $M$ extending the action of $\g g$.
\ssq We will be interested in the sequel by the following particular situation : the manifold $M$ is a vector space $V$, the action of $G$ is linear, and there is a classical r-matrix $r$ on $\g g$. We can as in the introduction view $J$ as a Lie algebra morphism from $\g g$ to the space of linear vector fields on $V$, and extend $J$ to a morphism $J^\bullet$ of Gerstenhaber algebras from $\Lambda(\g g)$ to $\wt\Lambda(V)$. In particular $J^2(r)$ defines a quadratic Poisson structure on $V$, and formula just above gives a quantization of this particular quadratic Poisson structure.
\paragraphe{III. Quadratic Poisson structures and  r-matrices}
\alinea{III.1. Some definitions}
We keep the notations of the introduction. The Gerstenhaber algebra $\widetilde \Lambda (V)$ can be written as~:
$$\widetilde \Lambda (V)=\bigoplus_{n\ge 0}\bigl(S^n(V)\otimes\Lambda^n(V)\bigr)[1]=\bigoplus_{n\ge 0}\widetilde \Lambda^n(V)[1].$$
A quadratic Poisson structure on $V$ can be defined as a bivector field $\Lambda$ in $\widetilde\Lambda^2(V)$ such that~:
$$[\Lambda,\, \Lambda]=0.$$
Let $\Lambda$ be an element of $\widetilde\Lambda^2(V)$, and let $r$ an element of $\Lambda^2(\g g)$ such that $J^2(r)=\Lambda$. It is then obvious that $[\Lambda,\Lambda]=0$ if and only if $J^3([r,r])=0$. If $n\ge 2$ then $J^2$ and $J^3$ have nontrivial kernels~: Precisely we have~:
$$\mop{dim}\mop{ker}J^2={n^2(n^2-1)\over 4}\espace{15}\hbox{and}\espace{15}
\mop{dim}\mop{ker}J^3={n^2(n^2-1)(5n^2-8)\over 36}.$$
\alinea{III.2. A counterexample in dimension 3}
\qquad
With the notations of \S\ I, an element of $\g g\wedge\g g$ can be written as~:
$$ r = \sum_{i,j,k,l=1}^{n} r_{i k}^{j l} E_{i j} \wedge E_{k l}. $$
Whe shall need for further calculations the following result~:
\prop {III.1}
Let $ r = \displaystyle{\sum_{i,j,k,l=1}^{n}} r_{i k}^{j l} E_{i j} \wedge E_{k l} $ be an
element of $ \g g \wedge \g g $ then $ [r,r] = 0 $ if and only if for any $i,j,k,l,m,p\in\{1,\ldots,n\}$ such that $(i,j)< (k,l) <(m,p)$ according to lexicographical order we have~:
$$ \sum_{d=1}^n r_{i k}^{d l} r_{d m}^{j p} - r_{m k}^{d l}r_{d i}^{p j} +
r_{k m}^{d p} r_{d i}^{l j}  - r_{i m}^{d p} r_{d k}^{j l} +  r_{m i}^{d j}
r_{d k}^{p l}  -  r_{k i}^{d j} r_{d  m}^{l p} =0 .$$
\dem
This proposition is a direct consequence of formula~:
$$ [ E_{ij} , E_{k l}] = 
\delta_{j k} E_{i l} - \delta_{ li} E_{ k j} $$ 
and the following lemma~:
\lemme {III.1}
Let $ \g h $ be a finite-dimensional Lie algebra and let $\uple XN$ be a basis of $ \g h $. 

If $ r = \displaystyle{\sum_{I,J=1}^N} r^{IJ} X_I \wedge
X_J $ is an element of $ \g h \wedge\g h $ ( $ r^{IJ} = - r^{JI} )$ then
$$ [r,r] = 4 \sum_{I,J,K,L=1}^N r^{IJ}r^{KL}[X_I,X_K] \wedge X_J \wedge X_L .$$
\ndem
\qed
{\sl Remark~:\/} We can directly show proposition III.1 using relation $(*)$ of beginning of \S\ II applied to elements of the dual basis of $\uple Xn$.
\medskip
\prop {III.2}
The Poisson structure $ \Lambda $ on $ \R^3 $ given by
$$ \Lambda = ( x_1^2 + \alpha x_2 x_3) \bord_2 \wedge \bord_3$$
with $\alpha \not = 0$
is not the image of a classical r-matrix by $ J^2 $.
\dem
An element $ r = {\displaystyle\sum_{i,j,k,l=1}^{n}} r_{i k}^{j l} E_{i j} \wedge E_{k l} $
of $ \g g \wedge \g g $ is parametrized by $36$ coefficients $r_{ik}^{jl}$. It has image $ \Lambda $ if and only if the following 18 equations are satisfied~:

 $$\eqalign{&r_{1 1 }^{1 2 } = r_{1 1}^{1 3} = r_{2 2}^{1 2} = r_{2 2}^{1 3} = 
r_{2 2}^{2 3} =  r_{3 3}^{1 2} = r_{3 3}^{1 3} = r_{3 3}^{2 3} = 0,\espace{15} 
r_{1 1}^{2 3} = 1 \cr
&r_{1 2}^{1 2} = r_{1 2}^{2 1},\espace{25} 
r_{1 2}^{1 3} = r_{1 2}^{3 1},\espace{25}
 r_{1 3}^{1 2} = r_{1 3}^{2 1}\cr
&  r_{1 3}^{1 3} = r_{1 3}^{3 1},\espace{25}
r_{1 2}^{2 3} = r_{1 2}^{3 2},\espace{25}
r_{1 3}^{2 3} = r_{1 3}^{3 2}\cr
&r_{2 3 }^{1 2} =r_{2 3}^{2 1},\espace{25} 
r_{2 3}^{1 3} = r_{2 3}^{3 1},\espace{25} 
r_{2 3}^{3 2} = r_{2 3}^{2 3} - \alpha \cr}$$
To lighten writing we rename the 18 remaining unknowns as follows~:
$$\eqalign{&r_{1 2}^{1 1} = a ,\esp r_{1 2}^{1 2} = b ,\esp  r_{1 2}^{1 3} = c ,\esp 
 r_{1 3}^{1 1} = d ,\esp  r_{1 3}^{1 2} = e ,\esp  r_{1 3}^{1 3} = f,\cr
 &r_{1 2}^{2 2} = g ,\esp r_{1 2}^{2 3} = h ,\esp r_{1 3}^{2 2} = i ,\esp 
 r_{1 3}^{2 3} = j ,\esp  r_{1 2}^{3 3} = k ,\esp  r_{1 3}^{3 3} =  l,\cr 
&r_{2 3}^{1 1} = m ,\esp  r_{2 3}^{1 2} = n,\esp r_{2 3}^{1 3} = p ,\esp 
 r_{2 3}^{2 2} = q ,\esp  r_{2 3}^{2 3} = r ,\esp  r_{2 3}^{3 3} = s.\cr}$$
The unknown $r$ must not be confused with the classical $r$-matrix on the whole. The context will not lead to any confusion. We have then~:
$$\eqalign{&r_{1 1}^{2 3} = 1 ,\cr
  &r_{1 2}^{2 1} = b ,\esp   r_{1 3}^{2 1} =  e,\esp  r_{1 2}^{3 1} = c ,\cr
&  r_{1 2}^{3 2} = h,\esp  r_{1 3}^{3 1} = f,\esp r_{1 3}^{3 2} = j ,\cr
&  r_{2 3}^{2 1} = n ,\esp  r_{2 3}^{3 1} = p,\esp r_{2 3}^{3 2} = r - \alpha,\cr}$$
and other $ r_{i k}^{j l} $ are equal to 0.
\ssq
If $r$ is such an element the equation $[r,r]=0$ develops according to Proposition III.2. into a system of 84 equations involving our 18 unknowns $a,b,\ldots, s$, given by the vanishing of the 84 coefficients of elements of the basis $E_{ij}\wedge E_{kl}\wedge E_{mn}$ of $\g g \wedge\g g\wedge\g g$. The 84 equations reduce to 66 thanks to the fact that we already have $[\Lambda, \Lambda]=0$. But we shall only consider $20$ of them, which will be sufficient for exhibiting the counterexample~:
\ssq
Let us order the $E_{ij}$'s lexicographically from first to $9$th, rename them accordingly ($A_1=E_{11}, A_2=E_{12},\ldots , A_9=E_{33}$), and labelize by $(x,y,z)$ the equation obtained by the vanishing of the coefficient of $A_x\wedge A_y\wedge A_z$. We shall consider precisely the following equations~:

$$\eqalign{& (1,2,5)\ ci-eh+n=0 \espace {70}(1,2,9)\ eh-ci+d=0 \cr
            & (1,3,5)\ cj-ek-a=0 \espace {70}(1,3,7)\ ce+f^2-ja-ld+m=0\cr
              &(1,3,9)\ ek-cj+p=0     \espace {70}(1,4,5)\ mh-nc=0\cr
        &(1,4,6)\ mk-pc=0       \espace {84}(1,5,6)\ \alpha c+nk-ph=0\cr
        &(1,7,8)\ en-im=0       \espace {84}(1,7,9)\ ep-jm=0   \cr
        &(1,8,9)\ \alpha e+ip-jn=0 \espace {64}(2,3,9)\ -3f+r+ik-jh=0\cr
   &(2,4,5)\ ag-b^2+nh-qc=0       \espace {42}(2,5,6)\ 2\alpha h+bh-rh+qk-cg=0\cr
  &(2,8,9)\ ej+ir+\alpha i-jq-fi=0\espace {23}(3,8,9)\ is+2\alpha j +el -fj-jr=0\cr
    &(4,5,7)\ pn-rm-bm+na=0 \espace {33}(5,6,9)\ -nk-rs+ hp+s(r-\alpha)=0\cr
        &(5,8,9)\ nj-ip+\alpha q=0 \espace {64}(6,8,9)\ ln-pj+qs-(r-\alpha)^2=0.\cr}$$

Consider the following two sums~:

$$\eqalign{(1,2,5) + (1,2,9): \esp &n+d=0\cr
(1,3,5) + (1,3,9): \esp &p-a=0.  \cr}$$
Hence $n=-d$ and $p=a$.
We will discuss the four cases $ a=d=0 $, $a=0$ and $d\not =0$, $a\not =0$ and
$d=0$, $ a\not =0 $ and $ d\not =0$.
\medskip
{\bf First case}~: $a=d=0$. Then looking successively at the following equations we get~:
\let\impl=\Longrightarrow

$$\eqalign{&(5,8,9) \impl        q=0
\espace {15}(6,8,9) \impl        r=\alpha
\espace {15}(1,5,6) \impl         c=0\cr
&(1,8,9)\impl          e=0
\espace{15}(2,4,5)\impl           b=0
\espace{15}(2,5,6)\impl           h=0\cr
&(4,5,7)\impl           m=0
\espace{12}(5,6,9)\impl           s=0
\espace{15}(1,3,7)\impl             f=0\cr
&(3,8,9)\impl           j=0
\espace{15}(2,8,9)  \impl         i=0
\espace{15}(2,3,9)  \impl         \alpha =0,\cr}$$
hence a contradiction to the hypothesis $\alpha\not=0$.
\smallskip
{\bf Second case}~: $a=0$ and $d\not= 0$ (hence $n\not =0$).
$$(1,4,6)\impl mk=0.$$
{\sl First subcase\/}~: $m=0$. Then~:
$$(1,4,5)\impl c=0\espace{15}
(1,7,8)\impl e=0\espace{15}
(1,2,5)\impl n=0,$$
hence a contradiction.
\smallskip
{\sl Second subcase\/}~: $m\not =0$, hence $k=0$.
$$(1,5,6)\impl c=0 \espace{15}
(1,4,5)\impl h=0\espace{15}
(1,2,5)\impl n=0,$$
hence a contradiction again.
\smallskip
{\bf Third case\/}~: $a\not =0$ and $d=0$ (hence $p\not =0$).
$$(1,4,5)\impl mh=0.$$
{\sl First subcase\/}~: $m=0$. Then~:
$$(1,4,6)\impl c=0\espace{15}
(1,7,9)\impl e=0\espace{15}
(1,3,5)\impl a=0,$$
hence a contradiction.
\smallskip
{\sl Second subcase\/}~: $m\not =0$, hence $h=0$.
$$(1,5,6)\impl c=0 \espace{15}
(1,4,6)\impl k=0\espace{15}
(1,3,5)\impl a=0,$$
hence a contradiction again.
\smallskip
{\bf Fourth case}~: $a\not =0$ and $d\not =0$.
\smallskip
{\sl First subcase\/}~: $m=0$.
$$(1,4,5)\impl c=0\espace{15}(1,7,8)\impl e=0\espace{15}(1,3,5)\impl a=0,$$
contradiction.
\smallskip
{\sl Second subcase\/}~: $m\not =0$.
$$(1,4,6)\impl k={a\over m}c\espace{15}(1,7,9)\impl j={a\over m}e
\espace{15} (1,3,5)\impl a=0,$$
contradiction.
\smallskip
This proves proposition III.2.
\qed
\alinea{III.3. Cartan-type quadratic Poisson structures}
\qquad Recall from \cite {DH} that the curl of a Poisson stucture $\Lambda=\sum_{i,j}\Lambda^{ij}\bord_i\wedge\bord_j$ is defined by~:
$$\mop{rot}\Lambda=\sum_{i,j}\bord_j\Lambda^{ij}.\bord_i.$$
It is a linear vector field (and hence can be viewed as an $n\times n$ matrix) when $\Lambda$ is quadratic. A quadratic Poisson structure is {\sl of Cartan type\/} if it can be written for some choice of coordinates as~:
$$\Lambda =\sum_{i,j=1}^nc_{ij}x_ix_j\bord_i\wedge\bord_j$$
with $c_{ji}=-c_{ij}$. J.P. Dufour and A. Haraki proved the following result~:
\th{III.3 \rm (Dufour - Haraki)}
Any quadratic Poisson structure the curl of which has eigenvalues $\lambda_i$ such that $\lambda_i+\lambda_j\not =\lambda_r+\lambda_s$ for any $(i,j,r,s)$ with $r\not =s$ and $\{i,j\}\not =\{r,s\}$ is of Cartan type.
\ndem
Such a Cartan-type Poisson structure is image by $J^{(2)}$, of a classical $r$-matrix, namely~:
$$r=\sum_{i,j}c_{ij}E_{ii}\wedge E_{jj}.$$

\paragraphe{References}
\bib{ABM}Arnal, D., Ben Amar, N., Masmoudi, M.: {Cohomology of good graphs and Kontsevich linear star products}. Lett. Math. phys. {\bf 48}, 291-306 (1999).
\bib{BR}Bhaskara, K.H., Rama, K.: {Quadratic Poisson structures}. J. Math. Phys. {\bf 32}, 2319-2322 (1991).
\bib{CFT}Cattaneo, A.S., Felder, G., Tomassini, L.: {From local to global deformation quantization of Poisson manifolds}. math.QA/0012228.
\bib{Di}Dito, G.: {Kontsevich star-product on the dual of a Lie algebra}. Lett. Math. Phys. {\bf 48}, 307-322 (1999).
\bib{Dr}Drinfel'd, V.G.: {On constant, quasiclassical solutions of the quantum Yang-Baxter equation}. Soviet Math. Dokl. {\bf 28}, No 3, 667-671 (1983).
\bib{DH}Dufour, J-P., Haraki, A.: {Rotationnels et structures de Poisson quadratiques}. C.R.Acad. Sci. {\bf 312}, 137-140 (1991).
\bib{ET}El Galiou, M.., Tihami, Q.: {Star-Product of a quadratic Poisson structure}. Tokyo J. Math. {\bf 19}, No 2, 475-498 (1996).
\bib{Gr}Grabowski, J.: {Abstract Jacobi and Poisson structures. Quantization and star-products}. J. Geom. Phys. {\bf 9}, 45-73 (1992).
\bib{K}Kontsevich, M.: {Deformation quantization of Poisson manifolds I}, math.QA/9709040.
\bib{Ka}Kathotia, V.: {Kontsevich's universal formula for deformation quantization and the Camp\-bell-Baker-Hausdorff formula}. math.QA/9811174.
\bib{LW}Lu, J-H., Weintein, A.: {Poisson Lie groups, dressing transformations and Bruhat decompositions}. J. Diff. Geom. {\bf 31}, 501-526 (1990).
\bib{LX}Liu, Z-J., Xu, P.: {On quadratic Poisson structures}. Lett. Math. Phys. {\bf 26}, 33-42 (1992).
\bib{OMY}Omori, H., Maeda, Y., Yoshioka, A.: {Deformation quantizations of Poisson algebras}. Contemp. Math. {\bf 179}, 213-240 (1994).
\bib{T}Takhtajan, L.A.: {Lectures on Quantum groups}. Nankai lect. notes in math. phys., 69-197 (1990).
\bib{V}Voronov, A.A.: {Homotopy Gerstenhaber algebras}. math.QA/9908040.
\bye